\documentclass[12pt,reqno]{amsart}

\usepackage[usenames]{color}
\usepackage{amssymb}
\usepackage{amsmath}
\usepackage{amsthm}
\usepackage{amsfonts}
\usepackage{amscd}
\usepackage{mathrsfs}
\usepackage{hyperref}
\usepackage{enumerate}
\usepackage{graphicx}
\DeclareMathOperator{\lcm}{lcm}

\def\stirlingtwo#1#2{S^{(2)}(#1,#2)}
\def\stirlingone#1#2{S^{(1)}(#1,#2)}

\def\sequencetwo{{S}^{(2)}}
\def\sequenceone{{S}^{(1)}}

\def\geq{\geqslant}

\theoremstyle{plain}

\newtheorem*{theorem*}{Theorem}
\newtheorem*{corollary*}{Corollary}
\newtheorem*{lemma*}{Lemma}
\newtheorem*{proposition*}{Proposition}

\theoremstyle{definition}

\newtheorem*{definition*}{Definition}
\newtheorem*{example*}{Example}
\newtheorem*{conjecture*}{Conjecture}

\theoremstyle{remark}

\newcommand\fix{{\rm{Fix}}}

\newcommand\divides\mid
\newcommand\smalldivides{\mathrel{\kern-2pt\kern3.5pt|}}
\newcommand\notdivides{\mathrel{\kern-3pt\not\!\kern4.3pt\bigm|}}
\newcommand\smallnotdivides{\mathrel{\kern-2pt\not\!\kern3.5pt|}}
\renewcommand{\le}{\leqslant}
\renewcommand{\ge}{\geqslant}

\begin{document}

\title[Patrick Moss]{Patrick Moss\\
25/10/1947--17/3/2024}

\author{Thomas Ward}
\address{Department of Mathematical Sciences,
Durham University, Durham, England}
\email{tbward@gmail.com}

\subjclass[2010]{Primary 37P35, 37C30, 11N32}

\date{\today}

\maketitle

\vspace{-6mm}

\begin{figure}[!h]
\begin{center}
\includegraphics[height=80pt]{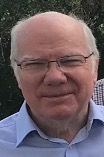}%
\end{center}
\end{figure}

There can be very few published mathematical
researchers with a~$39$ year gap in their
list of publications. Patrick (Pat) Moss, who completed a
doctorate with Graham Everest and me in~2003, is one of them.
He had a brief early period of activity
in ring theory
(see~\cite{MR0369424,MR0447309,MR0576182}).
After a long career in teaching he returned to
mathematics, where he had
a period of particularly
creative research in a form of abstract or
combinatorial dynamical systems
(see~\cite{MR4002553,pm,MR4394356}).

His contributions during this second period
of activity far exceeded what is visible
through his publications, in part because he
was consistently diffident about the depth of his
ideas.
While my focus here has been on the mathematical
contributions that Pat made, and his great
ingenuity, both Graham Everest and I liked him
very much and greatly enjoyed the all too short
period over which we were able to work closely with him.
I miss them both, as Graham himself passed
away in~2010~\cite{MR3105004}.

\section{Ring theory and Birkbeck College}

Pat grew up in Brentwood and attended school in
Chelmsford. He went on to study Mathematics
at the University of Surrey, graduating with
a BSc in~1968.
He went on to Birkbeck College as a postgraduate,
where he worked with Thomas Lenagan and
Stephen Ginn but it seems did not
complete his Masters degree.
In their first publication~\cite{MR0369424}
Moss and Ginn showed that a finitely generated finitely embedded projective module over a Noetherian ring is Artinian, answering a question of
Jategaonkar~\cite{MR0352170}.
They also showed that such a ring is Artinian if its right socle is either left or right essential.
In later work~\cite{MR0447309}
they used the notion of Artinian radical
and earlier work of Lenagan~\cite{MR0384862} to show that
a left and right Noetherian
ring with an Artinian classical quotient ring
is the direct sum of an Artinian
ring and a ring with zero socle.
Pat continued this line of enquiry,
trying to find necessary and/or
sufficient conditions for a Noetherian ring
to be an order in an Artinian ring,
with Lenagan. In joint work motivated in part
by earlier work with Ginn~\cite{MR0447309}
and work of Lenagan~\cite{MR0442008}
showing that several properties
of fully bounded Noetherian rings also hold
for Noetherian rings under the assumption that they have Krull dimension one,
they found simple conditions for the Artinian radical
of a Noetherian ring to split off
and under some more technical conditions on the rings found
criteria for the existence of an
Artinian quotient ring along the lines of the
results for Krull dimension one~\cite{MR0576182}.

\section{Palmers Sixth Form College 1970--2011}

Pat started teaching at Palmers Sixth Form College
in~1970, and this is the affiliation on his
joint work with Lenagan~\cite{MR0576182}.
He taught there for more than forty years,
inspiring many generations of students
with a love of Mathematics, and supporting
many students through the
Sixth Term Examination Paper (STEP) examinations
for access to universities that require it.
Pat met his wife Sally in~1994, and she survives him.

\section{Doctorate at UEA}

Pat completed a masters degree in Mathematics through
the Open University, graduating in~$2000$.
While studying there he met and befriended
Gerard McLaren, who was working
at British Telecom. Together they approached
the School of Mathematics at the University of East
Anglia to enquire about enrolling for PhDs part-time.
Gerry worked with Graham Everest and I on
properties of elliptic divisibility sequences,
leading to some results on effective bounds
for the appearance of primitive divisors
in certain families of such sequences~\cite{MR2220263}.
In the end Gerry decided not to continue with his
doctoral studies.

Pat also had Graham and me as supervisors,
and he looked at the notion of `realizability' for
integer sequences. This
purely combinatorial property had been studied
from different perspectives in multiple contexts,
but the starting point for Pat was some work
of Yash Puri~\cite{MR1873399}. A sequence denoted~$a=(a_n)$
of non-negative integers is said to be `realizable'
if there is a map~$T\colon X\to X$ for which
\[
a_n=\fix_T(n)=\vert\{x\in X\mid T^nx=x\}\vert
\]
for all~$n\ge1$. The sequence~$a$
is then said to be `realized' by the map~$T$.
Here~$X$ can be taken to simply be a countable
set with~$T$ a bijection or (via a suitable compactification) to be a compact
metric space with~$T$ a homeomorphism
or (less obviously)~$T$ can be taken to
be a~$C^{\infty}$ diffeomorphism and~$X$ an
annulus by work of Alastair Windsor~\cite{MR2422026}.

Part of the interest in this notion comes
from the congruence constraints that
realizability imposes---for
example, if~$a_1$ is odd then the fact that~$a_2-a_1$
counts the number of closed~$T$-orbits of length~$2$
forces~$a_2$ to also be odd.
In fact it is straightforward to see that the condition amounts
to the following: A sequence~$(a_n)$ is realizable if
and only if
\begin{itemize}
\item $\sum_{d\vert n}\mu\left(\frac{n}{d}\right)a_d\equiv0\pmod{n}$---the
`Dold' condition (D)
and
\item $\sum_{d\vert n}\mu\left(\frac{n}{d}\right)a_d\ge0$---the `Sign' condition (S)
\end{itemize}
for all~$n\ge1$, where~$\mu$ denotes the classical M{\"o}bius
function.
A recent overview of this property and settings in which
it arises may be found in the survey by
Byszewski, Graff, and the author~\cite{MR4332826}.

Yash Puri and~I had, for example, used the congruence~(D)
to study Fibonacci-like sequences, namely
those of the form
\[
(a_n)=(1,c,1+c,1+2c,2+3c,3+5c,\dots)
\]
satisfying the Fibonacci reccurence~$a_{n+2}=a_{n+1}+a_{n}$ for~$n\ge1$.
Our first result was that such a sequence
is realizable if and only if
the parameter~$c=3$,
meaning that the classical Lucas sequence is the
only solution of the Fibonacci recurrence that
is realizable~\cite{MR1866354}.
The argument uses well-known
congruences of the form~$F_{p-1}\equiv1$ modulo~$p$
for primes~$p\equiv\pm2$ modulo~$5$.
Much later Gregory Minton~\cite{MR3195758}
showed in essence that a linear recurrence sequence
satisfying the congruence~(D) must be a combination
of traces of powers of algebraic numbers, hugely
generalizing our observation.

Pat took these simple ideas and produced interesting new
insights in several different directions, many of them
pointing to further work.

\subsection{Local and Algebraic Realizability}

Pat looked at the~$p$-part of a sequence for each prime~$p$,
giving rise to the notion of `local' realizability: A
sequence is realizable locally at a prime~$p$ if the sequence
of~$p$-parts is itself a realizable sequence. It is
easy to see that a sequence that is locally realizable
at every prime is realizable, but Pat showed that the
converse also holds if the sequence is realized by
an endomorphism of a locally nilpotent group
(`nilpotently realizable').

This raised the question of understanding the
property of being realized by a group automorphism
more generally. If~$a$ is realized by a group
automorphism then in addition to the purely combinatorial
congruence conditions it must be a divisibility
sequence. Pat quickly found simple examples
to show that far more is needed:
Indeed, the sequence
\begin{equation}\label{cycle}
(1,1,1,1,6,1,1,1,1,6,\dots)
\end{equation}
which is realized by the permutation~$(12345)(6)$
on the set~$\{1,2,\dots,6\}$ is both a divisibility sequence
and a linear recurrence sequence, but it cannot
be realized by a group automorphism:
There is no group of order~$6$ with an automorphism
that cycles its non-identity elements.
What is clear about the sequence~\eqref{cycle}
is that both its~$2$-part and its~$3$-part are not realizable.
Pat's result~\cite[Th.~3.2.11]{pm} that nilpotent realizability
(realizability by an automorphism of a nilpotent
group) is
equivalent to everywhere locally nilpotent realizability
gives an explanation for this example, but the
problem of characterising algebraic realizability
combinatorially remains open, even for the restricted
class of linear recurrent divisibility sequences. These have
been classified by B\'{e}zivin, Peth\H{o}, and
van der Poorten~\cite{MR1081812} but which of them are algebraically
realizable is unclear.

\subsection{Realizability along Subsequences}

Pat developed an extraordinary insight into the
congruence conditions~(D) for realizability, and this
contributed greatly to
a series of conversations between Graham Everest,
Pat, and me about how sampling along a subsequence
of times might preserve or not preserve the
realizability property. Pat quickly established
a fundamental `powers' result: If~$(a_n)$ is realizable then~$(a_{n^k})$
is realizable for any~$k\ge1$~\cite[Th.~2.2.2]{pm}.
This raised the question of understanding what
kind of time-changing maps~$\mathbb{N}\to\mathbb{N}$ in place
of the map~$n\mapsto n^k$ also have
this property, but it was some years before any
of us were able to return to it.

In~2018 Sawian Jaidee visited me in Leeds, and we
realised that no polynomials apart from the
monomials already identified could have this
realizability along subsequences property. We
contacted Pat, and pieced together the following
way of thinking about these `time changes' that
preserve realizability---or, equivalently, act as
symmetries of the space of dynamical zeta functions.

\begin{definition*}[From Jaidee, Moss, and Ward~\cite{MR4002553}]\label{definitionofP}
For a dynamical system~$T\colon X\to X$ with~$\fix_{T}(n)<\infty$ for all~$n\geqslant1$,
define
\[
\mathscr{P}(X,T)=\{h\colon\mathbb{N}\to\mathbb{N}\mid\bigl(\fix_{T}(h(n))\bigr)
\mbox{ is a realizable sequence}\}
\]
to be the set of \emph{realizability-preserving time-changes for}~$(X,T)$.
Also define
\[
\mathscr{P}=\bigcap_{\{(X,T)\}}\mathscr{P}(X,T)
\]
to be the monoid of \emph{universally
realizability-preserving time-changes}, where
the intersection is taken over all systems~$(X,T)$
for which
\[
\fix_{T}(n)<\infty
\]
for all~$n\geqslant1$.
\end{definition*}

Pat had shown that the map~$n\mapsto n^k$ lies in~$\mathscr{P}$
for any~$k\in\mathbb{N}$, and our argument complemented this.
We also used the construction of certain
families of non-polynomial time
changes to show that the monoid~$\mathscr{P}$ has the
following properties (see~\cite{MR4002553}):
\begin{itemize}
\item the only polynomials in~$\mathscr{P}$ are the monomials, and
\item the monoid~$\mathscr{P}$ is nonetheless uncountable.
\end{itemize}
The first result says
that algebraically~$\mathscr{P}$ is very small;
the second that when viewed as a set~$\mathscr{P}$ is very large.

This circle of questions is still actively studied,
and recent work includes a complete description
of a set of (topological) generators for~$\mathscr{P}$
by Jaidee, Byszewski, and the author~\cite{jaidee2024generating}
and a geometric-combinatorial construction of the
map realizing~$(a_{n^2})$ by
Grzegorz Graff and Jacek Gulgowski (in progress).

\subsection{Polynomial Powers}

Pat also showed in~\cite[Cor.~2.1.6]{pm} that the terms of a realizable sequence
may themselves be raised to polynomial powers in the
following sense:
If~$h\in\mathbb{N}[x]$ and~$(a_n)$ is
realizable, then~$(a_n^{h(n)})$ is also
realizable. This is unexpected, and its
full ramifications have not been explored further.
What is expected is that polynomials are essentially
the only functions with this property.

\subsection{The Fibonacci Sequence}

Yash Puri and I had started our work by noticing
that the classical Fibonacci sequence
\[
(F_n)=(1,1,2,3,5,\dots)
\]
is not realizable. After the work
on the monoid~$\mathscr{P}$ with Sawian Jaidee,
Pat and I continued to correspond and he pointed out that
the sequence~$(5F_{n^2})$ is realizable.
This was an extraordinary
thing to notice, particularly as Pat did not use
computers to test ideas numerically, and of course
it is a very rapidly growing sequence.
The observation once again proved to be a
productive one.
First, this led naturally to the notion of
\emph{almost realizability}: A sequence~$(a_n)$
is almost realizable if there is a constant~$C$
such that~$(Ca_n)$ is realizable.
This notion became important later in work
of Miska~\cite{MR4361581} on the Stirling numbers.

\begin{definition*}[From Miska and Ward~\cite{MR4361581}]
Write~$\stirlingone{n}{k}$ for the Stirling
numbers of the first kind, defined for any~$n\ge1$ and~$0\le k\le n$
to be the number of permutations of~$\{1,\dots,n\}$
with exactly~$k$ cycles.
Write~$\stirlingtwo{n}{k}$ for
the Stirling numbers of the second kind,
defined for~$n\ge1$ and~$1\le k\le n$
to be the number of ways to partition
a set comprising~$n$ elements
into~$k$ non-empty subsets.
For each~$k\ge1$ define positive integer sequences
\begin{align*}
\sequenceone_k
&=
\bigl(\stirlingone{n+k-1}{k}\bigr)_{n\ge1}
\intertext{and}
\sequencetwo_k
&=
\bigl(\stirlingtwo{n+k-1}{k}\bigr)_{n\ge1}.
\end{align*}
\end{definition*}

\begin{theorem*}[From Miska and Ward~\cite{MR4361581}]
For~$k\ge1$ the sequence~$\sequenceone_k$ is
not almost realizable.
For~$k\le2$ the sequence~$\sequencetwo_k$ is realizable.
For~$k\ge3$ the sequence~$\sequencetwo_k$ is not
realizable, but is almost realizable
and the minimal constant multiplier~$C_k$ needed
is a divisor of~$(k-1)!$
for each~$k\ge1$.
\end{theorem*}

The multiplier needed for a given~$k\ge1$ in this
result is somewhat mysterious, as it is {\it{a priori}}
impossible to compute. It is by definition the
least common multiple of the denominators appearing
in an infinite sequence of rational numbers of
the form~$\frac1n\sum_{d\vert n}\mu\left(\frac{n}{d}\right)S^{(2)}(d+k-1,k)$.
If in calculating the first few terms
(or indeed the first few thousand terms) this
least common multiple reaches~$(k-1)!$ then this
certainly
is the minimal multiplier needed. If---as is usually
the case---this does not happen, then it is not clear
at any stage if the minimal value has really been
computed.
While some suggestions and speculations were
included in~\cite{MR4361581}, there is at this
stage very little knowledge about the
general term of the sequence of multipliers~$(C_k)$,
which begins
\[
(1, 1, 2, 6, 12, 60, 30, 210, 840, 2520, 1260, 13860, 13860, 180180,\dots).
\]

Going back to the Fibonacci sequence,
Pat and I were able to show that~$(F_{n^k})$
is almost realizable if~$k$ is even,
and is not almost realizable if~$k$ is odd.
Moreover, the constant needed for any even
power is~$5$, the discriminant of the Fibonacci
sequence~\cite{MR4394356}.

\begin{theorem*}[From Moss and Ward~\cite{MR4394356}]
If~$j$ is odd, then
the set of primes dividing denominators
of~$\frac{1}{n}\sum_{d\smalldivides n}\mu\bigl(\tfrac{n}{d}\bigr)F_{d^j}$
for~$n\in\mathbb{N}$ is infinite.
If~$j$ is even, then the sequence~$(F_{n^j})$ is
not realizable, but the sequence~$(5F_{n^j})$ is.
\end{theorem*}

The second natural direction of travel
starting with Pat's result about the
Fibonacci sequence sampled along the squares
is to ask how prevalent this phenomenon is,
and this has more or less been worked out in joint
work with Florian Luca~\cite{MR4590321}.
The details of this result are cumbersome,
but it is in essence a direct strengthening
of Pat's result on the Fibonacci sequence,
and the proof uses Binet's formula for the terms of a linear
recurrence sequence directly to replace the detailed
knowledge of congruence properties for the Fibonacci
sequence.

\begin{theorem*}[From Luca and Ward~\cite{MR4590321}]
Let~$u=(u_n)$ be a minimal linear recurrence sequence
satisfying
\[
u_{n+k}=a_1u_{n+k-1}+\cdots+a_k u_{n}
\]
with~$a_k\ne 0$,
and assume that the minimal
polynomial~$F$ of $u$ has only simple zeros.
Let~${\mathbb K}$ be the splitting field of~$F$
and~${\mathcal O}_{\mathbb K}$ be
its ring of integers. Let~$\Delta({\mathbb K})$ be the discriminant
of~${\mathbb K}$ and let~$\Delta(F)$ be the discriminant of~$F$. Let~$G$ be the Galois
group of~${\mathbb K}$ over~${\mathbb Q}$,
let~$e(G)$ be the exponent of~$G$, and let~$N$ be the order of~$G$.
\begin{itemize}
\item[{\rm{(i)}}] The sequence~$(Mu_{n^s})_{n\geq 1}$ satisfies
condition~{\rm{(D)}}
if~$M$ is a positive integer which is a multiple of~$\lcm(\Delta({\mathbb K}), \Delta(F))$ and~$s\geq N$ is a multiple of~$e(G)$.
\item[{\rm{(ii)}}]  Assume also that~$a_i\geq 0$ for~$i$
in~$\{1,\ldots,k\}$
and~$a_k\ne 0$, that
\[
(a_1,\ldots,a_k)\ne (0,0,\ldots,1),
\]
and that~$u_i\geq 1$ for all~$i$
in~$\{1,2,\ldots,k\}$.
Then the sequence
\[
(Mu_{n^s})_{n\geq 1}
\]
satisfies
condition~{\rm{(S)}}
whenever~$s=\ell e(G)$
where~$\ell\geq \ell_0$ is a sufficiently large
number which can be computed in terms of the sequence~$(u_n)_{n\geq 1}$.
\end{itemize}
\end{theorem*}

\subsection{Euler and Bernoulli Numbers}

Pat also showed realizability for several sequences
arising in arithmetic, two of which
seemed particularly interesting.
The `Euler numbers'~$(E_n)$
are defined by the formal relation
\[
\frac{2}{{\rm{e}}^t+{\rm{e}}^{-t}}=\sum_{n=0}^{\infty}E_n\frac{t^n}{n!}.
\]
By extending the classical Kummer congruences for
the Euler numbers, Pat was able to show the following,
recovering a result of Juan Arias de Reyna~\cite{MR2163516}.

\begin{theorem*}[From Moss~\cite{pm} {\&} Arias de Reyna~\cite{MR2163516}]
The sequence~$(\vert E_{2n}\vert)$ is realizable.
\end{theorem*}

The `Bernoulli numbers'~$(B_n)$
are similarly defined by the formal relation
\[
\frac{t}{{\rm{e}}^t-1}=\sum_{n=0}^{\infty}B_n\frac{t^n}{n!}.
\]
By making ingenious use of classical congruences
due to Adams, von Staudt-Clausen, and Kummer,
Pat was able to show something quite unexpected.
Write
\[
\left\vert
\frac{B_{2n}}{2n}
\right\vert
=
\frac{\tau_n}{\beta_n}
\]
with~$\gcd(\tau_n,\beta_n)=1$
and~$\tau_n,\beta_n\ge1$
for~$n\ge1$.

\begin{theorem*}[From Moss~\cite{pm}]
The sequence~$(\beta_n)$ is algebraically realizable.
\end{theorem*}

Pat also showed that the numerator sequence is
realizable (but cannot be algebraically realizable).

\begin{theorem*}[From Moss~\cite{pm}]
The sequence~$(\tau_n)$ is realizable.
\end{theorem*}

Indeed, he was able to show precisely how the
failure of algebraic realizability occurs.
Recall that a prime~$p$ is called irregular if
it divides the class number of the~$p$th
cyclotomic field~$\mathbb{Q}(\zeta_p)$ where~$\zeta_p$
denotes a primitive~$p$th root of unity.

\begin{theorem*}[From Moss~\cite{pm}]
The sequence~$(\tau_n)$ is not locally realizable
precisely at the set of irregular primes.
\end{theorem*}

\bibliographystyle{uwab}
\bibliography{refs}
\end{document}